\documentstyle[amscd,amssymb,verbatim,11pt]{amsart}
\newtheorem{Prop}{Proposition}[section]
\theoremstyle{definition}
\newtheorem{Def}[Prop]{Definition}
\newtheorem{Thm}[Prop]{\bf Theorem}
\newtheorem{Cor}[Prop]{Corollary}

\newtheorem{Lem}[Prop]{Lemma}
\theoremstyle{remark}
\newtheorem{Rem}[Prop]{Remark}
\numberwithin{equation}{section}

\def\Im{\operatorname{Im}}
\def\Tor{\operatorname{Tor}}

\def\dim{\operatorname{dim}}

\def\R{\text{\bf R}}

\def\Q{\text{\bf Q}}

\def\Z{\text{\bf Z}}

\def\N{\text{\bf N}}

\def\sL{\cal L}

\def\sP{\cal P}

\def\sS{\cal S}

\def\sC{\cal C}

\begin{document}

\title{Cohomological dimension of Markov compacta}
\author{A.N. Dranishnikov}
\address{Math Dept, University of Florida, Gainesville, FL 32611-8105, USA}
\email{dranish@@math.ufl.edu}
\thanks{The author was  supported in part by the NSF grant DMS0604494}



\date{}

\begin{abstract}
We rephrase Gromov's definition of Markov compacta, introduce a
subclass of Markov compacta defined by one building block and
study cohomological dimensions of these compacta. We show that for
a Markov compactum $X$, $\dim_{\Z_{(p)}}X=\dim_{\Q}X$ for all but
finitely many primes $p$ where $\Z_{(p)}$ is the localization of
$\Z$ at $p$. We construct Markov compacta of arbitrarily large
dimension having $\dim_{\Q}X=1$ as well as Markov compacta of
arbitrary large rational dimension with $\dim_{\Z_p}X=1$ for a
given $p$.
\end{abstract}
 \maketitle
\medskip
\medskip
\tableofcontents
\vfill
\pagebreak

\section{Introduction}

\subsection{Markov compacta}

Let $T$ be a rooted locally finite simplicial tree with the root
$x_0\in T$. For every vertex $x\in T$ by $T_x$ we denote the
subtree rooted at $x$, i.e. the tree with the vertices $y$ such
that the segment $[x_0,y]$ contains $x$. Gromov calls the tree $T$
{\it Markov}~\cite{Gr} if there are only finitely many (say, $k$)
isomorphism classes of rooted trees $T_x$. The name {\it Markov}
is given since the $k\times k$ transition matrix $M=(m_{ij})$
defines a Markov chain where $m_{ij}$ is the number of vertices of
the type $j$ neighboring the root  in a tree of type $i$.

A rooted tree can be viewed as the telescope of an inverse
sequence of finite spaces $S=\{K_i,\phi^{i+1}_i\}$ with $K_0=x_0$,
$|K_i|<\infty$. We call two points $x\in K_i$ and $y\in K_j$
equivalent if the inverse sequences
$$S_x=\{x\leftarrow(\phi^{i+1}_i)^{-1}(x)\leftarrow
(\phi^{i+2}_i)^{-1}(x)\leftarrow\dots\}$$
and
$$S_y=\{y\leftarrow(\phi^{j+1}_j)^{-1}(y)\leftarrow
(\phi^{j+2}_j)^{-1}(y)\leftarrow\dots\}$$
are isomorphic. Then Gromov's definition can be translated as
follows: An inverse sequence of finite spaces
$\{K_i,\phi^{i+1}_i\}$ is called {\it Markov} if this equivalence
relation on $\coprod_iK_i$ has finitely many classes.

This notion can be extend for sequences of higher dimensional
polyhedra.

\

\begin{Def}\label{markov spectra} Let $S=\{K_i,\phi^{i+1}_i\}$ be an inverse sequence of
simplicial complexes. We call this system {\it Markov} if the
following equivalence relation on the set of all simplices in
$\coprod K_i$ has finitely many classes: Two simplices
$\sigma\subset K_i$ and $\sigma'\subset K_j$ are equivalent if the
inverse sequences
$$S_{\sigma}=\{\sigma\leftarrow(\phi^{i+1}_i)^{-1}(\sigma)
\leftarrow(\phi^{i+2}_i)^{-1}(\sigma)\leftarrow\dots\}$$
and
$$S_{\sigma'}=\{\sigma'\leftarrow(\phi^{j+1}_j)^{-1}(\sigma')
\leftarrow(\phi^{j+2}_j)^{-1}(\sigma')\leftarrow\dots\}$$
are isomorphic.

A compactum $X$ is called {\it Markov} if it can be presented as
the limit of a Markov inverse system.
\end{Def}
\

Note that Markov inverse system consists of complexes of uniformly
bounded dimension. Hence Markov  compacta are always
finite-dimensional. Note that in a Markov system we have only
finitely many homeomorphism types of preimages of simplices
$(\phi^{i+1}_i)^{-1}(\sigma)$, we  call them {\it building blocks}.
Informally, Gromov defined Markov compacta as those which can be
built using finitely many building blocks. 
The Pontryagin surface $\Pi_2$ is a typical example of Markov compactum
which is constructed from
one building block. We recall the construction. Let
$f:M\to\Delta$ be a map of the M\"{o}bius band $M$ onto a
2-simplex which is the identity on the boundary. Fix a
sufficiently small triangulation on $M$. Take a triangulation of a
2-sphere and replace all its 2-simplices $\sigma$ by $M$ by means
of an identification $\partial\sigma\cong\partial M$. The
resulting space is supplied with natural projection onto $M$ glued
out of maps $f$. Then apply this procedure to the resulting space
and so on. We obtain an inverse system of polyhedra. The space
$\Pi_2$ is the limit space of this inverse sequence. Here the
building block is the M\"{o}bius band.

\subsection{One building block compacta}
Here we introduce a subclass of Markov compacta whose inverse
sequences can be obtained from one building block in some uniform
fashion.

A map $f:X\to Y$ is called {\it light} if the preimages of all
points $f^{-1}(y)$ are at most 0-dimensional. We call a simplicial
$n$-dimensional complex $K$ {\it a complex over an (oriented)
$n$-simplex $\Delta^n$} if there is a light simplicial map
$\chi:K\to\Delta^n$ (called a {\it characteristic map}). We denote
by $\beta K$ the barycentric subdivision of a simplicial complex
$K$. Note that $\beta K$  is a complex over $\Delta^n$ with the
characteristic map $\chi:\beta K\to\Delta^n$ defined on the
vertices of $\beta K$ as follows:
$\chi(b_{\sigma})=e_{\dim\sigma}$ where $b_{\sigma}$ denotes the
barycenter of a simplex $\sigma\subset K$ and $e_0,\dots,e_n$ are
the vertices of $\Delta^n$. The following Proposition is obvious.
\begin{Prop}\label{light}
Suppose that in the pull-back diagram
$$
\begin{CD}
K' @>>> L\\
@V\phi VV @VfVV\\
K @>\chi>> N\\
\end{CD}
$$
the map $f$ is simplicial and $\chi$ is light simplicial. Then
$K'$ is a simplicial complex and $\phi$ is simplicial map.
\end{Prop}

A triangulation $\tau$ of the simplex $\Delta^n$ is called {\it
symmetric} if it is invariant under the natural symmetric group
action on $\Delta^n$. Note that a symmetric triangulation $\tau$
on $\Delta^n$ induces a triangulation $\tau_K$ for every complex
$K$, $\chi:K\to\Delta^n$, over the simplex $\Delta^n$.

Let $f:L\to\Delta^n$ be a simplicial map of a finite complex $L$
onto the $n$-simplex $\Delta^n$ taken with a symmetric
triangulation $\tau$. Let $K_0$ be a complex over $n$-simplex with
the characteristic map $\chi_0:K_0\to\Delta^n$. By induction we
construct the following inverse sequence $\{K_i,\phi^{i+1}_i\}$ of
simplicial complexes over $\Delta^n$ with simplicial bonding maps
$\phi^{i+1}_i:K_{i+1}\to K_i$ With respect to some subdivision of
the triangulation on $K_i$.

Assume that $\chi_i:K_i\to\Delta^n$ is constructed. We define
$K_{i+1}$ as the pull-back of the diagram
$$
\begin{CD}
K_{i+1} @>\xi_{i+1}>> L\\
@V\phi^{i+1}_iVV @VfVV\\
K_i @>\chi_i>> \Delta^n.\\
\end{CD}
$$
The map $\chi$ is simplicial with respect the triangulation $\tau$
on $\Delta^n$ and the induced triangulation $\tau_{K_i}$. In view
of Proposition~\ref{light} $K_{i+1}$ is a simplicial complex and
the map $\phi^{i+1}_i$ is simplicial with respect to the
triangulation $\tau_{K_i}$ on $K_i$. We set the triangulation on
$K_{i+1}$ to be the first barycentric subdivision of $K_{i+1}$.
Then there is a natural characteristic map
$\chi_{i+1}:K_{i+1}\to\Delta^n$. The bonding map
$\phi^{i+1}_i:K_{i+1}\to K_i$ is simplicial with respect to
$\beta\tau_{K_i}$.

\begin{Def} The limit space $X$ of an
inverse sequence $\{K_i,\phi^{i+1}_i\}$ of complexes over the
$n$-simplex $\Delta^n$ defined above  is called a {\it  compactum
defined by the building block} $f:L\to\Delta^n$.
\end{Def}

\begin{Prop}
Every compactum defined by one building block is Markov.
\end{Prop}
\begin{pf}
Let $X$ be the limit space of an inverse sequence
$\{K_i,\phi^{i+1}_i\}$ of complexes over the $n$-simplex
$\Delta^n$ from the definition of compactum with one building
block and let  $f:L\to\Delta^n$ be the building block. We note
that simplices $\sigma_1\subset K_i$ and $\sigma_2\subset K_j$ are
equivalent (see Definition~\ref{markov spectra}) if
$\chi_i(\sigma_1)=\chi_j(\sigma_2)$ where $\chi_i:K_i\to\Delta^n$
and $\chi_j: K_j\to\Delta^n$ are the characteristic maps.
\end{pf}

In the case of the Pontryagin surface $\Pi_2$ we take $L$ to be
the M\"{o}bius band viewed as the mapping cylinder $M_g$ of a
2-fold covering map $g:S^1\to S^1$. We present the domain of $g$
as a 6-gon $S\simeq S^1$ and the range as a triangle $T\simeq
S^1$. Then we take $g$ simplicial. On the mapping cylinder of any
simplicial map always there is a triangulation on with no extra
vertices. We take such a triangulation on $L$ and define a
simplicial map $f:L\to\beta\Delta^2$ by an isomorphism taking $S$
onto $\beta(\partial\Delta^2)$ and by collapsing $T$ to the
barycenter $b_2$ of $\Delta^2$. Let $K_0$ be a 2-sphere with a
structure of a complex over the 2-simplex. Then $K_0$ and
$f:L\to\Delta^2$ define a compactum which is the Pontryagin
surface $\Pi_2$.

We recall that the Pontryagin surface $\Pi_2$ is 2-dimensional
with the rational dimension $\dim_{\Q}\Pi_2=1$. Mladen Bestvina asked me if
there are Markov compacta of dimension $n$ with rational dimension
one for arbitrary large $n$. Here we answer his question and give
an account of the cohomological dimension theory of Markov
compacta.

\subsection{Cohomological dimension}
Here is the summary of the cohomological dimension theory of
compacta (see \cite{Ku},\cite{Dr1},\cite{Dr2}). The cohomological
dimension of a space $X$ with coefficient group $G$ is defined as
follows:
$$
\dim_GX=\sup\{n|\ \check{H}^n(X,A;G)\ne 0\ \text{for some closed
subset}\ A\subset X\}.
$$
It is known for compact metric spaces that $\dim_GX\le n$ if and
only if the inclusion homomorphism
$\check{H}^n(X;G)\to\check{H}^n(A;G)$ is an epimorphism for every
closed subset $A\subset X$. The later is equivalent to the
condition that for every closed subset $A\subset $, every
continuous map $\phi:A\to K(G,n)$ to the Eilenberg-MacLane complex
has a continuous extension $\bar\phi:X\to K(G,n)$. By Bockstein
theorem to know the cohomological dimension of a compact space $X$
with respect to any abelian group it suffices to know it with
respect to the so-called {\it Bockstein groups} which are $\Q$,
$\Z_{(p)}$, $\Z_p$ and $\Z_{p^{\infty}}$ where $p$ runs over all
primes. Here $\Z_{(p)}$ is a localization of integers at $p$,
$\Z_p=\Z/p\Z$ and $\Z_{p^{\infty}}=\lim_{to}\Z_{p^k}$. In
particular, $\dim_{\Z}X=\sup\{\dim_{\Z_{(p)}}X\}$. The
cohomological dimension of compacta with respect to Bockstein's
groups is subject to restriction given by Bockstein Inequalities:
$$
\dim_{\Z_p}X-1\le\dim_{\Z_{p^{\infty}}}X\le\dim_{\Z_p}X;
$$
$$
\max\{\dim_{\Q}X,\dim_{\Z_p}X\}\le\dim_{\Z_{(p)}}X\le\max\{\dim_{\Q}X,\dim_{\Z_{p^{\infty}}}X+1\};
$$
$$
\dim_{\Z_{p^{\infty}}}X\le\max\{\dim_{\Q}X,\dim_{\Z_{(p)}}X-1\}.
$$
Let $\sigma$ denote the set of all Bockstein groups. There is
Realization Theorem \cite{Dr1},\cite{Dr2}, which states that for
every function $\beta:\sigma\to\N$ satisfying the Bockstein
inequalities there is a compact metric space $X$ with
$\dim_GX=\beta(G)$. A compactum $X$ is called {\it $p$-regular} if
$$\dim_{\Q}X=\dim_{\Z_{p^{\infty}}}X=\dim_{\Z_p}X=\dim_{\Z_{(p)}}X.$$
For a $p$-singular compactum $X$ the Bockstein inequalities split
into the following equality
$$\dim_{\Z_{(p)}}X=\max\{\dim_{\Q}X,\dim_{\Z_{p^{\infty}}}X+1\}$$
and the inequalities:
$$
\dim_{\Z_p}X-1\le\dim_{\Z_{p^{\infty}}}X\le\dim_{\Z_p}X.
$$
It is natural to suggest that Markov compacta could be
$p$-singular for only finitely many $p$. In this paper we prove
that $\dim_{\Z_{(p)}}X=\dim_{\Q}X$ for Markov compacta for all but
finitely many primes $p$.

 \

\subsection{Operations}
Clearly, the disjoint union and the product of two Markov compacta
are Markov. It is not the case for compacta defined by one block.
Nevertheless one can define these operation for such compacta.

Let $X$ and $Y$ be compacta defined by building blocks
$f:L\to\Delta^n$ and $g:N\to\Delta^m$, $n\ge m$, such that there
is an inclusion of the initial complexes $K_0\supset K_0'$. Let
$j:\Delta^m\to\Delta^n$ be the inclusion of the first $m$-face.
Then one can define their "sum" $X\# Y$ as the compactum generated
by the building block $f\cup (j\circ g):L\coprod N\to\Delta^n$.
This operation is most interesting when $n=m$ and
$K_0=K_0'=\Delta^n$.

The standard triangulation on the product of oriented simplices
$\Delta^n\times\Delta^m$ turns $\Delta^n\times\Delta^m$ into the
complex over $\Delta^{n+m}$. Let
$\chi:\Delta^n\times\Delta^m\to\Delta^{m+n}$ be its characteristic
map. Then we define the "product" $X\tilde\times Y$ of compacta
$X$ and $Y$ as the compactum defined by the building block
$\chi\circ(f\times g):L\times N\to\Delta^{n+m}$ and the initial
complex $K_0\times N_0$.

\subsection{Open problems}
Naturally, the compacta defined by one building block should have a fractal
structure. We recall that a compact set $F\subset\R^N$ is called
"self-similar" if there are finitely many similarities
$h_i:\R^N\to\R^N$ with the similarity coefficients $r_i<1$ such
that $F=\cup_ih_i(F)$ ~\cite{Ed}. 
Question: {\it Is every compactum defined by one
building block homeomorphic to self-similar subsets of $\R^N$?}

For compacta $X$ generated by one building block $f:L\to\Delta^n$ it would
be nice to obtain a formula for cohomological dimension or even
the formula for cohomology of $X$ in terms of $f$ in spirit of
those for Coxeter groups in terms of the nerve of Coxeter system
(\cite{Be},\cite{Dr3},\cite{D}).

\section{Cohomological dimension of one building block compacta}

\subsection{Restrictions on cohomological dimensions of Markov compacta}

Let $\phi:K\to K'$ be and let $A\subset N\subset K'$. We consider
the following condition:
$$
(*)^m_G\ \ \ \ \ \ im(f|_{f^{-1}(A)})^*\subset
im\{H^m(f^{-1}(N);G)\to H^m(f^{-1}(A);G)\}.
$$

We will use the notation $\phi\in(*)^m_G$ for saying that $\phi$
satisfies $(*)^m_G$ (for a certain pair $(N,A)$).

Easy diagram chasing yields the following.
\begin{Prop}\label{rel}
Let $\phi:K\to\Delta^{m}$ be a map to the $m$-simplex. Then
$\phi\in(*)^{m-1}_G$ for the pair $(\Delta^m,\partial\Delta^m)$ if
and only if the homomorphism
$\phi^*:H^m(\Delta^m,\partial\Delta^m;G)\to
H^m(K,\phi^{-1}(\partial\Delta^m);G)$ is nonzero.
\end{Prop}

\begin{Lem}\label{lemma1}
{\it Let $X=\lim_{\leftarrow}\{K_i;\phi^{i+1}_i\}$ be a Markov
compactum with $\dim_{G}X\le m$ for a principle ideal domain $G$.
Then for every $l\in\N\cup\{0\}$ there is $k$ such that the
inclusion $(*)^m_{G}$ holds for the map $\phi^{i+k}_i$ for all
$i\ge l$ for all pairs
$(\phi^i_{i-l})^{-1}(\sigma,\partial\sigma)$ where $\sigma$ is an
arbitrary simplex in $K_{i-l}$.}
\end{Lem}
\begin{pf}
Let $X$ be the limit of a Markov inverse sequence
$\{K_i,\phi^{i+1}_i\}$. Let $\sigma_i\subset K_i$ and
$\sigma_2\subset K_j$ be two equivalent in the sense of
Definition~\ref{markov spectra} simplices. Then we have
homeomorphic pairs
$$((\phi_i^{i+k})^{-1}(\sigma_1),(\phi_i^{i+k})^{-1}(\partial\sigma_1))
\ \ {\text{and}}\ \
((\phi_j^{j+k})^{-1}(\sigma_2),(\phi_j^{j+k})^{-1}(\partial\sigma_2))$$
for $k=0,1,\dots,\infty$. We take one representative
$\sigma\subset K_i$ for each equivalence class. By the definition
of Markov compactum there are only finitely many of them. Since
$A=(\phi^{i+1}_{i})^{-1}(\partial\sigma)$ is a finite complex, the
$G$-module $H^m(A;G)$ is finitely generated. Let $\{a_1,\dots,
a_s\}$ be a generating set.  Since $\dim_{G}X\le m$, the inclusion
$(\phi^{\infty}_{i+l})^{-1}(A)\subset
(\phi^{\infty}_{i+l})^{-1}(N)$ induces an epimorphism for
$m$-dimensional cohomology with coefficients in $G$ where
$N=(\phi^{i+1}_{i})^{-1}(\sigma)$. For every $j$ there is an
element $b_j\in\check{H}^m((\phi^{\infty}_{i+l})^{-1}(N);G)$ which
goes to $(\phi^{\infty}_{i+l})^*(a_j)$ under this inclusion
homomorphism. From the definition of \v{C}ech cohomology it
follows that there is $k_j$ such that
$(\phi^{i+l+k_j}_{i+l})^*(a_j)$ lies in the image of the
homomorphism induced by inclusion
$$(\phi^{i+l+k_j}_{i+l})^{-1}(A)\subset(\phi^{i+l+k_j}_{i+l})^{-1}(N).$$
We take $k$ greater than every $k_j$ for all equivalence classes.
\end{pf}
The converse to Lemma~\ref{lemma1} is true in the following form.
\begin{Lem}{}\label{lemma2}
{\it Suppose that a compact $X$ is presented as the inverse limit
of the sequence of $n$-dimensional polyhedra
$\{K_i,\phi^{i+1}_i\}$ supplied with triangulations $\tau_i$ such
that for every $j$
$$\lim_{i\to\infty}mesh(\phi^{j+i}_j(\tau_{j+i}))=0.$$ Assume that for every
$l$ there is $k$ such that the inclusion $(*)^m_G$ holds for the
maps $\phi^{i+k}_i$ for all $i\ge l$ for all pairs
$(\phi^i_{i-l})^{-1}(\sigma,\partial\sigma)$  where $\sigma$ is a
simplex in $K_{i-l}$. Then $\dim_GX\le m$.}
\end{Lem}
\begin{pf}
We show that given a continuous map $f:Y\to K(G,m)$ of a closed
subset $Y\subset X$ there is a continuous extension $\bar f:X\to
K(G,m)$. Since $K(G,m)$ is an ANE, there is $i_0$ and a map
$f':W\to K(G,m)$ of subcomplex $W\subset K_{i_0}$ which contains
$\phi^{\infty}_{i_0}(A)$ such that the composition
$f'\circ\phi^{\infty}_{i_0}|_A$ is homotopic to $f$. Here we used
the condition that the mesh of triangulations on $K_i$ tends to
zero. In view of the Homotopy Extension Theorem it suffices to
extend the map
$$g=f'\circ\phi^s_i|_{(\phi^s_i)^{-1}(W)}:(\phi^s_i)^{-1}(W)\to
K(G,m)$$ to $K_s$ for some $s$. Since $K(G,m)$ is
$(m-1)$-connected, there is an extension $f_m:W\cup
(K_{i_0})^{(m)}\to K(G,m)$. By induction on $i$ we define a number
$n_i$ and construct a map $$f_{m+i}:(\phi^{n_i}_{i_0})^{-1}(W\cup
(K_{i_0})^{m+i})\to K(G,m)$$ such that $n_{i}\ge n_{i-1}$ and
$f_{m+i}$ extends the map
$$f_{m+i-1}\circ\phi^{n_i}_{n_{i-1}}|_{(\phi^{n_i}_{i_0})^{-1}(W\cup
(K_{i_0})^{m+i-1})}.$$ Assume that $f_{m+i-1}$ is already
constructed. We take $k$ for $l=n_{i-1}-i_0$ from the condition of
Lemma and define $n_i=n_{i-1}+k$. For every $m+i$-dimensional
simplex $\sigma$ in $K_{i_0}\setminus W$ we consider the pair
$(\phi^{n_{i-1}}_{i_0})^{-1}(\sigma,\partial\sigma)$. By the
condition $(*)^m_G$ there is an extension
$\psi:(\phi^{n_{i}}_{i_0})^{-1}(\sigma)\to K(G,m)$ of the map
$f_{m+i-1}\circ\phi^{n_i}_{n_{i-1}}|_A$. The union of these
extensions for all $\sigma$ together with
$$f_{m+i-1}\circ\phi^{n_i}_{n_{i-1}}|_{(\phi^{n_i}_{i_0})^{-1}(W\cup
(K_{i_0})^{m+i-1})}$$ define $f_{m+i}$. Now the map $f_n$ is an
extension of the above map $g$ (for some $s$).
\end{pf}

\begin{Thm}\label{th1}
{\it For every Markov compactum $X$ there are only finitely many
primes $p_1,\dots, p_m$ such that
$\dim_{\Z_{(p_i)}}X\ne\dim_{\Q}X$.}
\end{Thm}
\begin{pf}
Let $X=\lim_{\leftarrow}\{K_i,\phi^{i+1}_i\}$ be a presentation of
$X$ from the definition of Markov compacta  and let
$\dim_{\Q}X=n$. Let $k=k(l)$ be from Lemma 2.1. Thus, the
condition $(*)_{\Q}^{n}$ holds for $\phi^{i+k}_i$ with
$(\phi^i_{i-l})^{-1}(\sigma,\partial\sigma)$ for all simplices
$\sigma$ in $K_{i-l}$ for all $i\ge l$. By the definition of
Markov compacta there are finitely many isomorphism types of
simplicial complexes in the family
$(\phi^{i+k}_{i-l})^{-1}(\sigma)$, $i\in \N\cup\{0\}$,
$\sigma\subset K_i$. Since all this complexes are finite, there is
$r_0$ such that for every prime $p> r_0$ the condition $(*)^n_G$
holds for $G=\Z_{(p)}$ for the pair
$(\phi^{i+k}_{i-l})^{-1}(\sigma,\partial\sigma)$ for every simplex
$\sigma$ in $K_{i-l}$ for all $i\ge l$. Lemma~\ref{lemma2} implies
the inequality $\dim_{\Z_{(p)}}X\le n$ for $p> r_0$. In view of
Bockstein inequality $\dim_{\Q}\le\dim_{\Z_{(p)}}$ we obtain
$\dim_{\Z_{(p)}}X=n$ for $p> r_0$.
\end{pf}

\subsection{Cohomological dimension of a complex over a simplex}

\begin{Def} Let $f:L\to\Delta^n$ be a map and let $G$ be an abelian
group. We define the {\it cohomological dimension} $cd_Gf$ of a
map $f$ with respect to the coefficient group $G$ to be the
minimal $m$ such that $f\in(*)^m_G$ for all pairs
$(\sigma,\partial\sigma)$ where $\sigma\subset\Delta^n$ is a
subsimplex. We define {\it the upper cohomological dimension}
$\overline{cd}_Gf$ of a map $f$ with respect to the coefficient
group $G$ to be the minimal $m$ such that the inclusion
homomorphism $H^m(f^{-1}(\sigma);G)\to
H^m(f^{-1}(\partial\sigma);G)$ is an epimorphism.
\end{Def}

Clearly, $cd_Gf\le\overline{cd}_Gf$.

Proposition~\ref{rel} implies the following.
\begin{Prop}\label{low bound}
Let $f:L\to\Delta^n$ be a map and let $G$ be an abelian group.
Then $cd_Gf$ is the maximal $k$ such that
$$f^*:H^k(\sigma^k,\partial\sigma^k;G)\to
H^k(f^{-1}(\sigma^k),f^{-1}(\partial\sigma^k);G)$$ is nonzero for
some $k$-face $\sigma^k\subset\Delta^n$.
\end{Prop}

\begin{Thm}\label{from above}
{\it Let $X$ be a compactum defined by the building block
$f:L\to\Delta^n$. Then $\dim_GX\le\overline{cd}_Gf$.}
\end{Thm}
\begin{pf}
The proof is similar to the proof of Lemma~\ref{lemma2}(1). Let
$\overline{cd}_Gf=m$. Given a continuous map $\psi:Y\to K(G,m)$ of
a closed subset $Y\subset X$, we construct a continuous extension.
We may assume that there is $i$ and a subcomplex $A_{i}\subset
K_{i}$ together with a map $g:A_{i}\to K(G,m)$ such that
$\phi^{\infty}_i(A)\subset A_{i}$ and $g\circ\phi^{\infty}_{i}|_A$
homotopic to $f$. Let
$$g':A_{i}\cup (K_{i})^{(m)}\to K(G,m)$$ be a continuous
extension. By the condition $\overline{cd}_Gf\le m$ we may assume
that for every $m+1$-simplex $\sigma$ in $K_{i}\setminus A_{i}$
there is an extension
$$g_{\sigma}^{m+1}:(\phi^{i+1}_{i})^{-1}(\sigma)\to K(G,m)$$ of the map
$$g'\circ\phi^{i+1}_i|_{(\phi^{i+1}_i)^{-1}(\partial\sigma)}:
(\phi^{i+1}_i)^{-1}(\partial\sigma)\to K(G,m).$$ The union of
$\cup_{\sigma}g_{\sigma}^{m+1}$ together with the composition
$g'\circ\phi^{i+1}_i)^{-1}|_{(\phi^{i+1}_i)^{-1}(A_i)}$ defines an
extension $ g^{m+1}:K_{i+1}^{(m+1}\cup(\phi^{i+1}_i)^{-1}(A_i)\to
K(G,m)$ of the map
$$g\circ\phi^{i+1}_i)^{-1}|_{(\phi^{i+1}_i)^{-1}(A_i)}.$$
Then for every $m+2$-simplex $\sigma$ in $K_i$ there is an
extension
$$g_{\sigma}^{m+2}:(\phi^{i+1}_{i})^{-1}(\sigma)\to K(G,m)$$ of the map
$g^{m+1}|_{(\phi^{i+1}_i)^{-1}(\partial\sigma)}$. The union of
$\cup_{\sigma}g_{\sigma}^{m+2}$ together with the map $g^{m+1}$
defines a continuous map
$g^{m+2}:K_{i+1}^{(m+2}\cup(\phi^{i+1}_i)^{-1}(A_i)\to K(G,m)$
extending $g\circ\phi^{i+1}_i)^{-1}|_{(\phi^{i+1}_i)^{-1}(A_i)}$
and so on. Repeating this procedure we will obtain a map
$g^{n}:K_{i+1}\to K(G,m)$ extending
$g\circ\phi^{i+1}_i)^{-1}|_{(\phi^{i+1}_i)^{-1}(A_i)}.$ Hence $
g^n\circ\phi^{\infty}_{i+1}$ is an extension of the map
$g\circ\phi^{\infty}_i|_A$. By the Homotopy Extension Theorem the
map $\psi$ has an extension.
\end{pf}
\begin{Rem}\label{better}
The argument of Theorem~\ref{from above} produces in fact a
slightly better inequality:
$$
\dim_GX\le\max\{cd_Gf,\overline{cd}_Gf-1\}.
$$
\end{Rem}

\subsection{Dimension over fields}

\begin{Lem}\label{L3}
{\it Suppose that a map $f:L\to\Delta^n$ simplicial with respect
to a symmetric triangulation $\tau$ on $\Delta^n$ induces an
epimorphism $$f_*:H_n(L,f^{-1}(\partial\Delta^n);R)\to
H_n(\Delta^n,\partial\Delta^n;R)$$ in the relative $n$-dimensional
homology with the coefficients in a ring $R$ with unit. Then for
every light simplicial map $\chi:K\to\Delta^n$ and any subcomplex
$A\subset K$ the induced homomorphism
$\phi_*:H_n(Z,\phi^{-1}(A);R)\to H_n(K,A;R)$ is an epimorphism
where $Z$ is the pull-back in the diagram
$$
\begin{CD}
Z @>\xi>> L\\
@V\phi VV @VfVV\\
K @>\chi>> \Delta^n.\\
\end{CD}
$$}
\end{Lem}
\begin{pf}
We define homomorphisms $f!$ of the simplicial chain complexes
with $R$-coefficients such that the diagram commutes (*):
$$
\begin{CD}
C_n(L) @>\partial >> C_{n-1}(L)\\
@Af!AA @Af!AA\\
C_n(\Delta^n)=R @>\partial>> C_{n-1}(\partial\Delta^n).\\
\end{CD}
$$
There is a relative $n$-cycle $z\in C_n(L)$ with the boundary
$\partial z\in C_{n-1}(f^{-1}(\partial\Delta^n))$ such that
$f_*([z])=1\in R=H_n(\Delta^n,\partial\Delta^n;R)$ where $[z]\in
H_n(L,f^{-1}(\partial\Delta^n);R)$ is the homology class of $z$.
Note that on the chain level we have a homomorphism $f_*:C_n(L)\to
C_n(\tau)$. We denote by $b: C_n(\Delta)\to C_n(\tau)$ the
subdivision homomorphism. Since $f_*([z])=[b(1\cdot\Delta^n)]$, by
dimensional reason there is only one element $b(1\cdot\Delta^n)$
in the relative homology class $[b(1\cdot\Delta^n)]$. Hence
$f_*(z)=b(1\cdot\Delta^n)$. We define $f!(\Delta)=z$. For an
$n-1$-face $\sigma\subset\Delta^n$ we define
$f!(\sigma)=pr_{\sigma}(\partial z)$ where
$pr_{\sigma}:C_{n-1}(L)\to C_{n-1}(f^{-1}(\sigma))$ the natural
projection. We define a homomorphisms $\phi!:C_*(K)\to C_*(Z)$ in
dimensions $n$ and $n-1$ on the generators $\sigma\subset K$ by
the formula
$$
\phi!(\sigma)=(\xi|_{\phi^{-1}(\sigma)})_*^{-1}f!\chi_*(\sigma).
$$
The following diagram is commutative (**):
$$
\begin{CD}
C_n(Z) @>\partial >> C_{n-1}(Z)\\
@A\phi!AA @A\phi!AA\\
C_n(K) @>\partial>> C_{n-1}(K).\\
\end{CD}
$$
It suffices to check the commutativity on the generators, i.e. the
equality $\partial \phi!(\sigma)=\phi!\partial(\sigma)$ for
$n$-simplices $\sigma\subset K$ taken with the coefficient $1\in
R$. The equality holds since for every $\sigma$ this diagram
contains a copy of the commutative diagram (*) with the
identification $\sigma=\Delta^n$.

Let $v\in C_n(K)$ be a relative cycle, i.e., $\partial v\in
C_{n-1}(A)$. From the commutativity of (**) it follows that
$\phi!(v)$ is a relative cycle with $\partial\phi!(v)\in
C_{n-1}(\phi^{-1}(A))$. We note that $\phi_*(\phi!(v))=b_K(v)$
where $b_K:C_n(K)\to C_n(\tau_K)$ is the subdivision homomorphism
and $\tau_K$ is the triangulation on $K$ induced from $\tau$ by
means of the map $\chi$. Then $\phi_*([\phi!(v)])=[b_K(v)]=[v]$
for the relative homology classes.
\end{pf}

\begin{Lem}\label{L4}
{\it Let $X=\lim_{\leftarrow}\{K_i,\phi^{i+1}_i\}$ be a compactum
defined by a building block $f:L\to\Delta^n$.

(1) Suppose the inequality $cd_Gf<n$ holds. Then $\dim_GX<n$.

(2) Let $F$ be aa additive group of a field and let $cd_Ff=n$.
Then $\dim_FX=n$.}
\end{Lem}
\begin{pf}
We may assume that $K_0=\Delta^n$.

(1) Follows from Remark~\ref{better}.

(2) In view of Proposition~\ref{rel} the homomorphism
$$f^*:H^n(\Delta^n,\partial\Delta^n;F)\to H^n(L,f^{-1}(L);F)$$ is
nontrivial. Since $F$ is a field, the dual homomorphism
$$f_*:H_n(L,f^{-1}(\partial\Delta^n);F)\to
H_n(\Delta^n,\partial\Delta^n;F)=F$$ is nontrivial and hence, it
is an epimorphism. Denote by $\partial
X=(\phi^{\infty}_0)^{-1}(\partial\Delta^n)$.  By induction using
Lemma~\ref{L3} we can construct a sequence $$v_i\in
H_n(K_i,(\phi^{i}_0)^{-1}(\partial\Delta^n);F)$$ such that
$(\phi^i_{i-1})_*(v_i)=v_{i-1}$ and $v_0=1\in F=
H_n(\Delta^n,\partial\Delta^n;F)$. Thus, we construct a nontrivial
$n$-dimensional relative \v{C}ech $F$-homology class on
$(X,\partial X)$. This implies that $H_{n}(X,\partial X;F)\ne 0$
for the Steenrod homology.  By the Universal Coefficient Theorem
over a field we obtain $\check{H}^n(X,\partial X;F)\ne 0$ for the
\v{C}ech cohomology with $F$-coefficients. Hence $\dim_FX\ge n$,
which contradicts to the assumption.
\end{pf}

We recall that for compact spaces $X$ there are two possibilities
for the dimension of the $n$-th power:
$$ \dim X^n=n\dim X \ \text{for all}\ n\ \ \text{or}\ \ \dim
X^n=(n-1)\dim X+1.
$$
We conjecture that all Markov compacta are of the first type. In
the support of the conjecture we present the following.

\begin{Thm}\label{bolt}
{\it Suppose that an $n$-dimensional compactum $X$ is defined by a
building block $f:L\to\Delta^n$ has dimension $n$. Then $\dim
X^k=kn$ for all $k$.}
\end{Thm}
\begin{pf}
Let  $p$ be a prime such that $\dim_{\Z_{(p)}}X=n$. In view of
Lemma~\ref{L4}(1) we have $cd_{\Z_{(p)}}f=n$. This implies that
$$f^*:H^n(\Delta^n,\partial\Delta^n;\Z_{(p)})\to
H^n(L,f^{-1}(\partial\Delta^n);\Z_{(p)})$$ is nontrivial. If
$\dim_{\Q}X=n$, the Theorem follows form the Kunneth formula over
the field. So, we assume that $\dim_{\Q}X<n$. Then by
Lemma~\ref{L4}(2) the inequality $cd_{\Q}f<n$ holds.  Hence
$$f^*:H^n(\Delta^n,\partial\Delta^n;\Q)\to
H^n(L,f^{-1}(\partial\Delta^n);\Q)$$ is a zero homomorphism. Hence
the image of $f^*$ with $\Z_{(p)}$-coefficient is a $p$-torsion
group. By the Universal Coefficient Formula for $\Z_p$ as a module
over $\Z_{(p)}$ we obtain that the homomorphism
$$f^*:H^n(\Delta^n,\partial\Delta^n;\Z_{p})\to
H^n(L,f^{-1}(\partial\Delta^n);\Z_{p})$$ is nontrivial. Hence
$cd_{\Z_p}f=n$.  By Lemma~\ref{L4}(2) $\dim_{\Z_p}X=n$. Then
$kn=k\dim X\ge\dim X^k\ge\dim_{\Z_p}X^k=k\dim_{\Z_p}X= kn$.
\end{pf}

\

\subsection{Symmetric building blocks}

The group of all permutations on $n$ elements is denoted by $S_n$.
There is a natural action of $S_{n+1}$ on the $n$-simplex
$\Delta^n$. A compactum defined by a building block
$f:L\to\Delta^n$ is called {\it symmetric} if there is an action
on $L$ of the symmetric group $S_{n+1}$ and the map $f$ is
$S_{n+1}$-equivariant.
\begin{Thm}\label{coh dim}
{\it Let $X$ be a symmetric compactum with a building block
$f:L\to\Delta^n$. Then for every field $F$ there are the
inequalities
$$cd_Ff\le\dim_FX\le\overline{cd}_Ff.$$}
\end{Thm}
\begin{pf}
In view of Theorem~\ref{from above} it suffices to prove only the
first inequality. Let $cd_Ff=m$. By Proposition~\ref{low bound}
there is an $m$-face $\sigma\subset\Delta^n$ such that
$cd_Ff|_{f^{-1}(\sigma)}=m$. Since $f$ is symmetric, we may assume
that $\sigma$ is the first $m$-face $\Delta^m\subset\Delta^n$. Let
$Y$ be a compactum defined by the building block
$f|_{f^{-1}(\sigma)}:f^{-1}(\sigma)\to\Delta^m$. We claim that
there is an embedding $Y\subset X$. Let $\{K_i,\phi^{i+1}_i\}$ and
$\{N_i,\psi^{i+1}_i\}$ be inverse systems for $X$ and $Y$ from the
definition compacta generated by one building block. Without loss
of generality we may assume that $K_0=\Delta^n$ and $N_0=\Delta^m$
with the identities as the characteristic maps. By induction we
construct an embedding of inverse sequences $N_i\subset K_i$. The
imbedding $\Delta^m\subset\Delta^n$ induces an imbedding
$N_1\subset K_1$. Since $\Delta^m$ is the first face and the
characteristic maps on $K_1$ and $N_1$ are defined by means of the
barycentric subdivision and the ordering of vertices of $\Delta^n$
and $\Delta^m$, we have that the restriction $\chi_1|_{N_1}$ is
the characteristic map for $N_1$. Therefore there is an embedding
of $N_2\subset K_2$ defined by the pull-back diagram from the
definition of compacta defined by one building block, and so on.

By Lemma~\ref{L4}(2) $\dim_FY=m$. Hence $\dim_FX\ge m$.
\end{pf}

\begin{Def} Let $f:L\to\Delta^n$. A {\it symmetrization} of $f$ is a map
$\tilde f:L\times S_{n+1}\to\Delta^n$ defined by the formula:
$\tilde f(x,s)=s(f(x))$.
\end{Def}

It is easy to see that the map $\tilde f$ is $S_{n+1}$-equivariant
with respect to the action on $L\times S_{n+1}$ generated by
multiplication in $S_{n+1}$ from the left and with natural action
on $\Delta^n$. The following is obvious.
\begin{Prop}\label{dim of symm}
$cd_Gf=cd_G\tilde f$ and $\overline{cd}_Gf=\overline{cd}_G\tilde
f$.
\end{Prop}

We note that the compactum $\tilde X$ obtained from the
symmetrization $\tilde f$ of $f:L\to\Delta^n$ is homeomorphic to
the sum
$$\#_{a\in S_{n+1}}X_a$$ where $X_a$ is generated by $a\circ f$.

\begin{Prop}\label{symmetrization}
For every compactum $X$ defined by a building block
$f:L\to\Delta^n$ there is a symmetric compactum $\tilde X$ defined
by the building block $\tilde f:L\times S_{n+1}\to\Delta^n$ that
contains $X$ as a subspace.
\end{Prop}
\begin{pf}
The embedding $X\subset\tilde X$ is induced by the diagram
$$
\begin{CD}
L @>{x\mapsto (x,e)}>> L\times S_{n+1}\\
@VfVV @V\tilde fVV\\
\Delta^n @>=>> \Delta^n\\
\end{CD}
$$
where $e$ is the unit in $S_{n+1}$.
\end{pf}

\section{Markov compacta with low rational dimension}

The main results of this section is the following theorem:
\begin{Thm}\label{main 1}
{\it For every $n\in\N$ and $k\le n$, for every finite set of
primes ${\sL}$ there is a (symmetric) compactum $X$ defined by one
building block $f_n:L_n\to\Delta^n$ with dimensions $\dim X=n$ and
$\dim_{\Z[\frac{1}{p}]}X=k$ for $p\in{\sL}$ for every $k\le n$.}
\end{Thm}

\subsection{Rational dimension $\ge$ 2}

First we prove this theorem for $k>1$.

Let $K_0 @>g_0>> K_1 @>g_1>> K_2 @>g_2>> \dots$ be a direct
sequence. The telescope $T(\{g_i\})$ generated by this sequence is
the quotient space $\coprod M_{g_i}/\sim$ where $M_{g_i}$ is the
mapping cylinder of the map $g_i:K_i\to K_{i+1}$ and the
equivalence relation $\sim$ identifies $K_i\subset M_{g_i}$ with
$K_i\subset M_{g_{i-1}}$.

Let ${\sS}$ be a subset of the set ${\sP}$ of all prime numbers.
The standard construction of a localization $X_{({\sS})}$ of a
space $X$ at ${\sS}$ uses the Postnikov tower. Sullivan's original
construction of the localization for a simply connected CW
complexes [Su] defines $X_{({\sS})}$ as an infinite telescope
$T(\{\nu_i\})$ of the direct sequence of simply connected
complexes
$$
K_0 @>\nu_0>> K_1 @>\nu_1>> K_2 @>\nu_2>>\dots
$$
with $K_0=X$, the localization map $l:X\to T(\{\nu_i\})$ equal to
the inclusion, and $\dim K_i=\dim X$ for all $i$. In this case we
say that an $n$-dimensional space $X$ {\it admits a localization
by means of a direct sequence of $n$-dimensional polyhedra}. Thus,
the Sullivan's construction gives such a localization for every
simply connected complex.
\begin{Prop}\label{factor}
Let $K$ be a finite simply connected simplicial complex of $\dim
K=n$ and let $p\in{\sP}$. Then there exists a finite simply
connected $n$-dimensional simplicial complex $K'$ and a map
$g:K\to K'$, simplicial with respect to some iterated barycentric
subdivision of $K$, such that the localization map $l:K\to
K_{({\sP}\setminus\{p\})}$ is homotopically factored through $f$,
$l\sim\xi\circ f$ and
$$
g_*:H_*(K;\Z_p)\to H_*(K';\Z_p)
$$
is zero homomorphism.
\end{Prop}
\begin{pf}
According to the above we may assume that
$K_{({\sP}\setminus\{p\})}=T(\{\nu_i\})$ for a sequence of simply
connected $n$-dimensional simplicial complexes
$$
K_0 @>\nu_0>> K_1 @>\nu_2>> K_2 @>\nu_3>>\dots
$$
with $K_0=K$ and $\nu_i$ simplicial with respect to some iterated
barycentric subdivision of $K_i$. First we note that
$H_*(T(\{\nu_i\});\Z_p)=0$. Since
$$T(\{\nu_i\}_{i=0}^{\infty})=\lim_{\to}T(\{\nu_i\}_{i=0}^j),$$ for
every element $\alpha\in H_*(K;\Z_p)$ there is $j(\alpha)$ such
that the image of $\alpha$ is zero in the finite telescope
$T(\{\nu_i\}_{i=0}^{j(\alpha)})$. Since $K$ is a finite complex,
there is $m$ such that the inclusion $$l_m:K\to
T(\{\nu_i\}_{i=0}^{m})$$ induces zero homomorphism for the mod $p$
homology. Note tat the telescope $T(\{\nu_i\}_{i=0}^{j(\alpha)})$
can be deformed to the space $K_m$. Let
$r:T(\{\nu_i\}_{i=0}^{j(\alpha)})\to K_m$ be the resulting
retraction. We take $K'=K_m$ and $g=r\circ l_m$. Note that $g$ is
a simplicial map for $s$-iterated barycentric subdivision of $K$
for sufficiently large $s$.
\end{pf}
Given a map $g:X\to Y$ we denote by $M_g$ and $C_g$ the mapping
cylinder and the mapping cone respectively. By $\Sigma X$ we
denote the suspension over $X$ and by $CX$ the cone over $X$.
\begin{Prop}\label{Prop mod p iso}
Let $g:K\to K'$ be as in Proposition~\ref{factor}, and let
$q:C_g\to C_g/K'=\Sigma K'$ be the projection. Then
$$
q_*:H_{n+1}(C_g;\Z_p)\to H_{n+1}(\Sigma K';\Z_p)
$$
is an isomorphism.
\end{Prop}
\begin{pf}
Consider the diagram generated by the exact sequence of homology
with coefficients in $\Z_p$ and the inclusions
$(CK,K)\to(C_g,M_f)\leftarrow (C_g,K')$.
$$
\CD
0=H_{n+1}(K') @>>> H_{n+1}(C_g) @>q_*>> H_{n+1}(C_g,K') @>>> H_n(K')\\
@ V=VV  @V=VV  @ V=VV  @ V=VV\\
0=H_{n+1}(M_g)  @>>> H_{n+1}(C_g) @>q_*>> H_{n+1}(C_g,M_g)
@>\partial>> H_n(M_g)\\
@AAA @AAA  @A=AA  @ Ag_*AA\\
\dots @>>>  0  @>>> H_{n+1}(CK,K)  @>\partial'>> H_n(K)\\
\endCD
$$

By Proposition~\ref{factor}, $g_*=0$. Hence $\partial=0$ and the
result follows. \
\end{pf}

Let us fix a prime $p$ and a natural number $k>1$. We define a
collection of building blocks $\{f_n:L_n\to\Delta^n\}$, $n\ge k$
by induction on $n$ such that each complex $L_n$ is simply
connected $n$-dimensional. We define the $L_k=\Delta^k$ and
$f_k=id_{\Delta^k}$. Assume that simply connected $i$-dimensional
simplicial complexes $L_i$ together with $f_i:L_i\to\Delta^i$ are
defined for $i< n$ such that $\dim L_i=i$ and the maps $f_i$ are
simplicial with respect to a symmetric triangulation $\tau^i$ of
$\Delta^i$. Let $\chi_n:\beta(\partial\Delta^n)\to\Delta^{n-1}$ be
the characteristic map. Denote by $\widetilde{\partial\Delta^n}$
the pull-back of the diagram
$$
\begin{CD}
\widetilde{\partial\Delta^n} @>>>L_{n-1}\\
@Vf_{n-1}'VV @Vf_{n-1}VV\\
\partial\Delta^n @>\chi>>\Delta^{n-1}.\\
\end{CD}
$$
Note that the space $\widetilde{\partial\Delta^n}$ is simply
connected and $n-1$-dimensional.  By Proposition~\ref{light} the
map $f_{n-1}'$ is simplicial with respect to the triangulation
$\tau_{n-1}'$ on $\partial\Delta^n$ induced from $\tau^{n-1}$ by
means the map $\chi$. Let $g:\widetilde{\partial\Delta^n}\to K'$
be a map from Proposition~\ref{factor} for the complex
$K=\widetilde{\partial\Delta^n}$. It is simplicial with respect to
the $s$-iterated barycentric subdivision $\beta^sK$ of $K$ for
some $s$. We define $L_n$ as the mapping cylinder $M_{g}$ where
$g$ is taken from Proposition~\ref{factor} for the complex
$\widetilde{\partial\Delta^n}$. We define the triangulation on
$\Delta^n$ as the cone $\tau^n=cone(\beta^s\tau_{n-1}')$ of the
$s$-iterated barycentric subdivision of the triangulation
$\tau_{n-1}'$ of the boundary $\partial\delta^n$. Note that it is
symmetric. The mapping cylinder of a simplicial map admits a
triangulation which coincides with the triangulation $\beta^sK$ on
$\widetilde{\partial\Delta^n}$. We fix such triangulation on $L_n$
and define $f_n:L_n\to\Delta^n$ as the simplicial map with respect
to $\tau^n$ that takes all vertices from $L_n$ which are not in
$\widetilde{\partial\Delta^n}$ to the cone vertex and coincides
with $\beta^sf_{n-1}'$ on $\widetilde{\partial\Delta^n}$.  Note
that the complex $L_n$ is simply connected and $n$-dimensional.

\begin{Lem}\label{mod p iso}
For every $n$-dimensional simplicial complex $K$ over $\Delta^n$,
$\chi:K\to\Delta^n$, the projection $\pi:\widetilde K\to K$ in the
pull-back diagram
$$
\begin{CD}
\widetilde{K} @>>>L_{n}\\
@V\pi VV @Vf_{n}VV\\
K @>\chi>>\Delta^{n}\\
\end{CD}
$$
induces an isomorphism $\pi_*:H_n(\widetilde K;\Z_p)\to
H_n(K;\Z_p)$.
\end{Lem}
\begin{pf}
We prove it by induction on $n$.

Consider the diagram generated by the mod $p$ homology and the
mapping $\pi:(\widetilde K,\pi^{-1}(K^{(n-1)}))\to (K,K^{(n-1)})$:

$$
\CD
0 @>>>H_n(\widetilde K) @>>> H_n(\widetilde
K,\pi^{-1}(K^{(n-1)}))
@ >>> H_{n-1}(\pi^{-1}(K^{(n-1)}))\\
@. @V\pi_*VV @V{\alpha}VV @ V{\beta}VV\\
0 @>>> H_n(K) @>>> H_n(K,K^{(n-1)})@ >>> H_{n-1}(K^{(n-1)}).\\
\endCD
$$
By the construction we can identify $\pi^{-1}(K^{(n-1)})$ with the
pull back of the diagram
$$
\CD
\widetilde{K^{(n-1)}} @>>>L_{n-1}\\
@V\pi|_{...}VV @Vf_{n-1}VV\\
K^{(n-1)} @>\chi_n\circ\chi>> \Delta^{n-1}.\\
\endCD
$$
By induction assumption $\beta$ is an isomorphism. We show that
$\alpha$ is an isomorphism as well and apply the Five Lemma. We
note that $\alpha$ is induced by the map $\bar\pi:\widetilde
K/\widetilde{K^{(n-1)}}\to K/K^{(n-1)}$ which is the wedge of maps
$\bar f_{n}:\widetilde{\Delta^n}/\widetilde{\partial\Delta^n}
\to\Delta^n/
\partial\Delta^n$ induced by $f_n$. By the construction
$\widetilde{\Delta^n}/\widetilde{\partial\Delta^n}= C_g$ where
$g:\widetilde{\partial\Delta^n}\to K'$ is from
Proposition~\ref{factor} and the projection $\bar f_{n}$ can be
factored as
$$
C_{g} @>q>> \Sigma\widetilde{\partial\Delta^n} @>\Sigma
(f_{n-1}')>> \Sigma\partial\Delta^n.
$$
By induction assumption $\Sigma f_{n-1}'$ induces isomorphism of
$n$-dimensional homology. Then Proposition~\ref{Prop mod p iso}
implies that $\bar f_{n}$ induces an isomorphism.
\end{pf}
\begin{Cor}\label{cd mod p}
{\it For every $n$, $cd_{\Z_p}f_n=n$.}
\end{Cor}
\begin{pf}
Apply Lemma~\ref{mod p iso} to the diagram generated by the map
$f_n:(L_n,f^{-1}_n(\partial\Delta^n)\to(\Delta^n,\partial\Delta^n)$
to obtain that $f_n$ induces nontrivial homomorphism of relative
cohomology with coefficients in $\Z_p$. Hence, $cd_{\Z_p}f_n\ge
n$.
\end{pf}
\begin{Prop}\label{top 1/p epi}
For every simplex $\Delta^n$ the inclusion
$\widetilde{\partial\Delta}\subset L_n$ induces an epimorphism of
k-dimensional cohomology with coefficients in $\Z[\frac{1}{p}]$.
\end{Prop}
\begin{pf}
This follows from the fact that the inclusion of a space $X$ to
its localization $X_{{\sP}\setminus\{p\}}$ induces  an isomorphism
for cohomology with $\Z[\frac{1}{p}]$ coefficients. Since the
localization map for $\widetilde{\partial\Delta}$ is homotopy
factored through the inclusion
$\widetilde{\partial\Delta}\subset\widetilde\Delta$ the required
statement follows.
\end{pf}
\begin{Lem}\label{1/p epi}
{\it For every $n$-dimensional simplicial complex $K$ over
$\Delta^n$, $\chi:K\to\Delta^n$ and every subcomplex $N\subset K$,
the inclusion $\widetilde N\subset\widetilde K$ induces an
epimorphism of $k$-dimensional cohomology with coefficients in
$\Z[\frac{1}{p}]$ where $\widetilde K$ is the pull-back in the
diagram
$$
\begin{CD}
\widetilde{K} @>>>L_{n}\\
@V\pi VV @Vf_{n}VV\\
K @>\chi>>\Delta^{n}\\
\end{CD}
$$
and $\widetilde N=\pi^{-1}(N)$.}
\end{Lem}
\begin{pf}
By induction on $n$. Lemma is true for $n=k$ by the dimensional
reason. Let $\phi:\widetilde N\to K(\Z[\frac{1}{p}],k)$ be a map.
We construct an extension $\bar\phi:\widetilde K\to
K(\Z[\frac{1}{p}],k)$. We note that $\pi^{-1}(K^{n-1})$ is the
pull-back in the diagram
$$
\begin{CD}
\pi^{-1}(K^{(n-1)}) @>>>L_{n-1}\\
@V\pi|_{...} VV @Vf_{n-1}VV\\
K^{(n-1)} @>\chi_n\circ\chi|_{K^{(n-1)}}>>\Delta^{n-1}.\\
\end{CD}
$$
By induction assumption there is an extension
$\phi':\pi^{-1}(K^{n-1})\to K(\Z[\frac{1}{p}],k)$ of the map
$$\phi|_{\pi^{-1}(N^{(n-1)})}:\pi^{-1}(N^{(n-1)})\to
K(\Z[\frac{1}{p}],k).$$ For every $n$-simplex $\sigma\subset K$ the
pair $(\pi^{-1}(\sigma),\pi^{-1}(\partial\sigma))$ is homeomorphic
to the pair $(L_n,f_n^{-1}(\Delta^n))$. By Proposition~\ref{top
1/p epi} there is an extension $$\phi_{\sigma}:\pi^{-1}(\sigma)\to
K(\Z[\frac{1}{p}],k)$$ of $\phi'$ restricted to
$\pi^{-1}(\partial\sigma)$. The union of $\phi_{\sigma}$ for
$\sigma\subset K\setminus Int(N)$ together with $\phi$ gives us a
required extension $\bar\phi$.
\end{pf}

The above construction can be summarized in the following.

\begin{Lem}\label{block1}
{\it Let $p$ be a prime. Then for every $n\in\N$ and $k$ with
$2\le k\le n$ there are an $n$-dimensional simplicial complex
$L_n$ and a map $f_n:L_n\to\Delta^n$ simplicial for some symmetric
triangulation of $\Delta^n$ such that $cd_{\Z_p}f_n=n$,
$cd_{\Q}f=k$ and $\overline{cd}_{\Z[\frac{1}{p}]}f\le k$.

Moreover, if ${\sL}=\{p_1,\dots,p_s\}$ is a finite set of primes,
then for the above $k$ and $n$ there is a map $f_n:L_n\to\Delta^n$
simplicial for some symmetric triangulation of $\Delta^n$ such
that $\ cd_{\Z_{p_1\dots p_s}}f_n=n$ , $cd_{\Q}f=k$ and
$\overline{cd}_{\Z[\frac{1}{p}]}f\le k$ for all $p\in{\sL}$.}
\end{Lem}
\begin{pf}
The case of one $p$ is presented above.

In the general case we replace $p$ by the product $p_1\dots p_s$.
The construction and the proof remain the same.
\end{pf}

PROOF OF THEOREM~\ref{main 1} $(k>1)$. Let
${\sL}=\{p_1,\dots,p_k\}$. By passing to the symmetrization, we
may assume that the map $f_n:L_n\to\Delta^n$ in Lemma~\ref{block1}
is symmetric (see Proposition~\ref{dim of symm}). Let $X_n$ denote
a compactum defined by $f_n:L_n\to\Delta^n$. By Theorem~\ref{coh
dim} $\dim_{\Z_p}X_n\ge n$, $p\in{\sL}$, $\dim_{\Q}X_n\ge k$ and
$\dim_{\Z[\frac{1}{\sL}]}X_n\le k$. The first inequality implies
that $\dim X_n=n$. The other two inequalities together with the
Bockstein inequalities imply $\dim_{\Z[\frac{1}{\sL}]}X_n=k$. \qed

\subsection{Rational dimension one}

The following changes are needed to run the construction for the
Theorem~\ref{main 1} for $k=1$. First in the presence of the
fundamental group the localization does not necessarily exists. So
Proposition~\ref{factor} must be changed. Still there is a
localization for homology, i.e. a map $X\to \bar X$ such that
$H_*(X)\to H_*(\bar X)$ is the localization homomorphism. The
problem here is that this localization is not necessarily given by
the direct system of complexes of the same dimension ($= \dim X$).
To make Proposition~\ref{factor} working we map our complex to a
complex of this type by a map that induced an epimorphism in
1-dimensional cohomology with the localized coefficient group.

\begin{Prop}\label{hom localization}
Let $p$ be a prime number. Let $L$ denote a finite product
$T^m\times\prod_{i=1}^sK(G_i,1)$ of Eilenberg-MacLane complexes
where $T^m$ is the $m$-torus, $G_i=\Z_{q_i^{m_i}}$ where $q_i$ is
prime and $K(G_i,1)$ is a complex with finite skeletons in all
dimensions for all $i$ . Then for every $n$ the $n$-skeleton
$L^{(n)}$ admits a homology localization at ${\sP}\setminus\{p\}$
by means of a direct system of $n$-dimensional polyhedra.
\end{Prop}
\begin{pf}
Let $p:S^1\to S^1$ be a map of degree $p$ and let $p^m:T^m\to T^m$
be the product of $m$ copies of $p$. We define
$\gamma_i:K(G_i,1)\to K(G_i,1)$ as follows. If $q_i\ne p$ we
define $\gamma_i=id$, if $q_i=p$ we take $\gamma_i$ to be a map to
a vertex in $K(G_i,1)$. Consider the map
$$\gamma=p^m\times\prod_{i=1}^s\gamma_i:
T^m\times\prod_{i=1}^sK(G_i,1)\to
T^m\times\prod_{i=1}^sK(G_i,1).$$ Clearly, $\gamma(L^{(n)})\subset
L^{(n)}$. It is easy to check that the iteration of
$\gamma|_{L^{(n)}}$ localizes the free part of the homology and
the torsion part.
\end{pf}
\begin{Prop}\label{new map}
Let $K$ be a  finite simplicial complex of $\dim K=n>1$. Then
there is a map $\psi:K\to K_0$ to an $n$-dimensional complex $K_0$
such that $$\psi^*:H^1(K_0;\Z[\frac{1}{p}])\to
H^1(K;\Z[\frac{1}{p}])$$ is an epimorphism and $K_0$ admits a
homology localization at ${\sP}\setminus\{p\}$ by means of a
direct sequence of finite $n$-dimensional polyhedra.
\end{Prop}
\begin{pf}
We attach finitely many 2-cells to $K$ to make the fundamental
group abelian. Let $N$ denote a new complex and let $j:K\to N$ be
the inclusion. Clearly, $\dim N=n$. There is a map $\alpha:N\to L$
where $L$ is as in Proposition~\ref{hom localization} that induces
an isomorphism of the fundamental groups. By the Universal
Coefficient Theorem $\alpha$ induces an isomorphism of
1-cohomology with coefficients in $\Z[\frac{1}{p}]$. We may assume
that $\alpha$ lands in $L^{(n)}$. Now take $K_0=L^{(n)}$ and
$\psi=\alpha\circ j$.\
\end{pf}
The following is a modification of Proposition~\ref{factor}.
\begin{Prop}\label{new factor}
Let $K$ be a finite simplicial complex of $\dim K=n$ and let $p\in\sP$.
Then there exists a finite $n$-dimensional simplicial complex
$K'$ and a map $g:K\to K'$, simplicial with respect to some
iterated barycentric subdivision
of $K$,
such that
$$
(1)\ \ \ \ g_*:H_*(K;\Z_p)\to H_*(K';\Z_p)
$$
is zero homomorphism and
$$
 (2)\ \ \ \ g_*:H^1(K';\Z[\frac{1}{p}])\to H^1(K;Z[\frac{1}{p}])
$$
is an epimorphism.
\end{Prop}
\begin{pf}
Take $\psi:K\to K_0$ from Proposition~\ref{new map} and consider a direct system
$$
K_0 @>\nu_0>> K_1 @>\nu_2>> K_2 @>\nu_3>>\dots
$$
that localizes homology of $K_0$. Then
$\lim_{\to}H_*(K_i;\Z_p)=0$. Take $i$ such that
$(\nu^0_i)_*:H_*(K_0;\Z_p)\to H_*(K_i;\Z_p)$ is zero homomorphism.
Then take $K'=K_i$ and $g=\nu^0_i\circ\psi$. Then (1) holds.

The homomorphism $(\nu^0_i)_*$ with coefficients in $\Z[\frac{1}{p}]$
is a monomorphism as a left devisor of the localization isomorphism.
Since $Ext$ term is zero in the Universal Coefficient Theorem over the ring
$\Z[\frac{1}{p}]$ for 1-dimensional cohomology, we have that
$(\nu^0_i)^*$ is an epimorphism for 1-dimensional cohomology with coefficients
in $\Z[\frac{1}{p}]$.

We may assume that $g$ is simplicial with respect to some iterated barycentric
subdivision of $K$.
\end{pf}
The prove of the following proposition differs from the proof of
Proposition~\ref{Prop mod p iso} only by the reference to
Proposition~\ref{new factor} instead of Proposition~\ref{factor}.
\begin{Prop} Let $g:K\to K'$ be as in Proposition~\ref{new factor}, and let
$q:C_g\to C_g/K'=\Sigma K'$ be the projection. Then
$$
q_*:H_{n+1}(C_g;\Z_p)\to H_{n+1}(\Sigma K;\Z_p)
$$
is an isomorphism.
\end{Prop}

We have constructed a map $g$ such that Lemma~\ref{block1} holds
for $k=1$ with the same proof. The proof of Theorem~\ref{main 1}
for $k=1$ goes without changes.

\section{Markov compacta with low mod $p$ dimension}

\begin{Thm}\label{main 2}
For every $n\in\N$ and $k\le n$, for every finite set of primes
${\sL}$ there is an $n$-dimensional compactum $X$ defined by a
building block $f_n:L_n\to\Delta^n$ with $\dim_{\Z_p}X=k$ for
$p\in{\sL}$ and every $k\le n$.
\end{Thm}

\subsection{Mod $p$ dimension $\ge 2$}

We denote by ${\sC}$ the class (Serre class) of torsion abelian
groups~\cite{Sp}.

\begin{Lem}\label{block2}
{\it Let $p$ be a prime. Then for every $n,k\in\N$ with $2\le k\le
n$ there are an $n$-dimensional simplicial complex $L_n$ and a map
$f_n:L_n\to\Delta^n$ simplicial with respect to some symmetric
triangulation of $\Delta^n$ such that
$cd_{\Z_p}f_n=\overline{cd}_{\Z_p}f_n= k$ and $cd_{\Q}f_n= n$.

Furthermore, for every finite set of primes
${\sL}=\{p_1,\dots,p_s\}$ for every $n,k\in\N$ with $2\le k\le n$
there are an $n$-dimensional simplicial complex $L_n$ and a map
$f_n:L_n\to\Delta^n$ simplicial with respect to some symmetric
triangulation of $\Delta^n$ such that
$cd_{\Z_p}f_n=\overline{cd}_{\Z_p}f_n= k$  for every $p\in {\sL}$
and $cd_{\Q}f_n= n$.}
\end{Lem}
\begin{pf}
In the case when $s>1$ we set $p=p_1\dots p_s$. Let $k\ge 2$ be fixed.
We use induction on $n$. We additionally assume the following:

(1) $L_n$ is simply connected,

(2) for every complex $K$ over the
simplex $\Delta^n$,
$\chi:K\to\Delta^n$ the homomorphism $f'_*:H_*(\widetilde K;\Q)\to H_*(K;\Q)$ induced by
$f'$ from the pull-back diagram
$$
\CD
\widetilde K @>\chi'>> L_n\\
@Vf'VV @Vf_nVV\\
K @>\chi>> \Delta^n.\\
\endCD
$$
is an isomorphism,

(3) for every complex $K$ over the
simplex $\Delta^n$,
$\chi:K\to\Delta^n$, and every subcomplex $A\subset K$,
the inclusion homomorphism $j^*:H^k(\widetilde K;\Z_p)\to
H^k(f^{-1}(A);\Z_p)$ is an epimorphism.

If $n=k$, we take $f_n=id_{\Delta^n}$.
Then all conditions are satisfied.

Now assume that $f_n:L_n\to\Delta^n$ is constructed and $n>k$.
Let $\chi:\partial\Delta^{n+1}\to\Delta^n$ be a characteristic map
of the barycentric subdivision of $\partial\Delta^{n+1}$.
Let $\widetilde{\partial\Delta^{n+1}}$ be the pull-back in the diagram
$$
\CD
\widetilde{\partial\Delta^{n+1}} @>\chi'>> L_n\\
@Vf'VV @Vf_nVV\\
\partial\Delta^{n+1} @>\chi>> \Delta^n.\\
\endCD
$$
By induction assumption $f'_*:H_*(\widetilde{\partial\Delta^{n+1}};\Q)\to
H_*(\partial\Delta^{n+1};\Q)$ is an isomorphism.
Hence with $\Z$ coefficients the induced homomorphism $f'_*$ is a
${\sC}$-isomorphism. Note that
$\widetilde\partial\Delta^{n+1}$ is 1-connected.
Then by the Mod ${\sC}$ Hurewicz Theorem
$f'_{\#}:\pi_n(\widetilde\partial\Delta^{n+1})\to\pi_n(\partial^{n+1})=\Z$
is a ${\sC}$-isomorphism. Hence
$\pi_n(\widetilde\partial\Delta^{n+1})=\Z\oplus Tor$.
Let $g:S^n\to\widetilde\partial\Delta^{n+1}$ be a map that represents
an element $p\in\Z$.
We define $L_{n+1}$ as the mapping cone of $g$. There is a natural map
$f_{n+1}:L_{n+1}\to\Delta^{n+1}$ which coincides with $f'$ over
$\partial\Delta^{n+1}$.

First we show that $\overline{cd}_{\Z_p}f_{n+1}\le k$. Since $f_n$
is the identity above the $k$-skeleton of $\Delta^{n+1}$, we have
that $\overline{cd}_{\Z_p}f_{n+1}\ge k$. When $n=k$, the space
$L_n$ is the mapping cone $C_g$ of a map $g:S^k\to S^k$ of degree
$p$ and the inclusion $\partial\Delta^{n+1}\cong S^k\subset C_g$
induces a monomorphism of $k$-homologies with $\Z_p$ coefficients.
Hence it induces an epimorphism of $k$-cohomologies with $\Z_p$
coefficients which implies that $\overline{cd}_{\Z_p}f_{n+1}\le
k$. Let $n>k$ and let $\sigma\subset\Delta^{n+1}$ be a face. If
$\sigma\subset\partial\Delta^{n+1}$, then the restriction of
$\chi$ to $\sigma$ makes it to be a complex over the simplex
$\Delta^n$. Then by the induction assumption (3),
$H^k((f')^{-1}(\sigma);\Z_p)\to
H^k((f')^{-1}(\partial\sigma);\Z_p)$ is an epimorphism. If
$\sigma=\Delta^{n+1}$, the inclusion homomorphism
$$H^k(L_{n+1};\Z_p)\to
H^k(\widetilde{\partial\Delta^{n+1}};\Z_p)$$ is an epimorphism,
since $L_{n+1}$ is obtained from
$\widetilde{\partial\Delta^{n+1}}$ by attaching one cell of
dimension $\ge k+1$ and hence the $k+1$-dimensional skeleton of
$L_{n+1}$ equals the $k+1$-dimensional skeleton of
$\widetilde{\partial\Delta^{n+1}}$. Finally, we note that
$cd_{\Z_p}f_{n+1}\ge k$. Hence
$cd_{\Z_p}f_n=\overline{cd}_{\Z_p}f_{n+1}=k$.

By construction, the inclusion
$\widetilde{\partial\Delta^{n+1}}\subset L_{n+1}$ induces zero
homomorphism in $n$-dimensional $\Q$-homology and hence in
$\Q$-cohomology. This implies that $cd_{\Q}f_{n+1}= n+1$.

Next we verify the conditions (1)-(3).

(1). Using induction assumption it is not difficult to show that
$\widetilde{\partial\Delta^{n+1}}$ is simply connected. Then
$L_{n+1}$ is simply connected by construction.

(2). Let $\nu:K\to\Delta^{n+1}$ be a light simplicial map. Then by
the construction of $L_{n+1}$, the
restriction
$$\nu'=\chi\circ\nu|_{K^{(n)}}:K^{(n)}\to\Delta^n$$ is a
light simplicial map such that $(f')^{-1}(K^{(n)})$ is the
pull-back of
$$
K^{(n)} @>\nu'>>\Delta^n @<f_n<< L_n.
$$
By the induction assumption
$$
(f'|_{...})_*:H_*((f')^{-1}(K^{(n)});\Q)\to H_*(K^{(n)};\Q)$$ is
an isomorphism. Consider the diagram generated by the exact
sequences of pairs and the map $f':(\widetilde
K,(f')^{-1}(K^{(n)}))\to(K,K^{(n)})$.
$$ \CD H_i(\widetilde K;\Q)
@>>>H_i(\widetilde K, (f')^{-1}(K^{(n)});\Q) @>>>
H_{i-1}((f')^{-1}(K^{(n)});\Q)\\
@Vf'_*VV @V{\psi}VV @V(f'|_{...})_*VV\\
H_i(K;\Q) @>>> H_i(K,K^{(n)};\Q) @>>> H_{i-1}(K^{(n)};\Q).\\
\endCD
$$
By the construction
$$\xi=(f_{n+1})_*:H_*(L_{n+1},\widetilde{\partial\Delta^{n+1}};\Q)\to
H_*(\Delta^{n+1},\partial\Delta^{n+1};\Q)$$ is an isomorphism.
Therefore, $\psi$ is an isomorphism as the direct sum of $\xi$. By
Five Lemma $f'_*$ is an isomorphism.

(3). It suffices to show that every map $\phi:(f')^{-1}(A)\to
K(\Z_p,k)$ has an extension $\bar\phi:\tilde K\to K(\Z_p,k)$. As
in the proof of (2) we may use the induction assumption to
construct an extension $\phi':\widetilde{K^{(n)}}\to K(\Z_p,k)$ of
the map $\phi|_{(f')^{-1}(A\cap K^{(n)})}$. Since the inclusion
$\widetilde{\partial\Delta^{n+1}}\subset L_{n+1}$ induces an
epimorphism in $k$-dimensional mod $p$ cohomology, there is an
extension $\bar\phi:\tilde K\to K(\Z_p,k)$ of the map
$\phi\cup\phi':(f')^{-1}(A)\cup\widetilde{K^{(n)}}\to K(\Z_p,k)$.
\end{pf}

PROOF OF THEOREM~\ref{main 2}($k>1$). Let $X_n$ be a (symmetric)
compactum defined by the block $f_n:L_n\to\Delta^n$. By
Theorem~\ref{coh dim}  implies that $\dim_{\Z_p}X_n=k$ and
$\dim_{\Q}X_n\ge n$. Hence $\dim X_n=n$. In case of $p=p_1\dots
p_s$ we have $\dim_{\Z_{p_i}}X_n\le\dim_{\Z_p}X_n\le k$. Since the
inequalities $cd_{G}f_n\ge k$ for any $G$, we obtain that
$\dim_{\Z_{p_i}}X_n=k$ for all $i$.\qed

\subsection{Mod $p$ dimension one}

To prove Lemma~\ref{block2} for $k=1$ we need a sequence of
results.

Let ${\sC}_p$ denote the Serre class of $p$-torsion groups.

\begin{Prop}\label{modification}
Let $X_0\subset X_1\subset\dots\subset X_n$ be a sequence of cell complexes
such that $X_0$ is finite and each $X_{i+1}$ is obtained from $X_i$ by attaching
finitely many (possibly no) $(i+1)$-dimensional cells. Suppose that
$\Tor H_n(X_0)\in {\sC}_p$ and $\dim X_0\le n$. Then $\Tor H_n(X_n)\in{\sC}_p$.
\end{Prop}
\begin{pf}
Since $\dim(X_{n-1}/X_0)\le n-1$, the exact sequence of the pair
$(X_{n-1}, X_0)$ implies that $H_n(X_0)=H_n(X_{n-1})$. Consider
the exact sequence of the pair $(X_n,X_{n-1})$:
$$
0=H_{n+1}(X_n,X_{n-1})\to H_n(X_{n-1}) @>i>> H_n(X_n)\to H_n(X_n,X_{n-1})\to.
$$
Note that $H_n(X_n,X_{n-1})=H_n(\vee S^n)=\oplus\Z$. Then
$H_n(X_n)\subset \Im  i\cong (\oplus\Z)\oplus \Tor H_n(X_0)$.
Thus, $\Tor H_n(X_n)\subset\Tor H_n(X_0)$.
\end{pf}

Let $q:X\to Y$ be the projection onto the orbit space of
$G$-action for a finite group $G$. We denote by $\tau:H_*(Y)\to
H^*(X)$ the homology transfer. Note that if the $G$-action is
free, then $q_*\tau_*$ is multiplication by $|G|$.

\begin{Prop}\label{Hurewicz}
Let $X$ be a complex with $\pi_1(X)=\oplus_{i=1}^{m}\Z_p$,
and $\pi_i(X)=0$ for $2\le i\le n$. Then $H_i(X)\in{\sC}_p$ for $i\le n$.
\end{Prop}
\begin{pf}
Let $q:\bar X\to X$ be the universal cover. Then $q_*\tau_*$ is the
multiplication
by $p^m$. By the Hurewicz theorem $H_i(\bar X)=0$ for $i\le n$. Thus the
homomorphism of multiplication by $p^m$ in $H_i(X)$ is zero. It means that
$H_i(X)\in{\sC}_p$.
\end{pf}
\begin{Prop}\label{futher modif}
Let $X$ be an $n$-dimensional compact polyhedron such that
$\pi_1(X)=\oplus_{i=1}^{m}\Z_p$, $\pi_i(X)=0$ for $2\le i\le n$
and $\Tor H_n(X)\in{\sC}_p$. Then by attaching finitely many $n+1$
cells to $X$, it is possible to construct a complex $Y$ such that
$H_i(Y)\in{\sC}_p$ for all $i$.
\end{Prop}
\begin{pf}
We chose a basis $a_1,\dots, a_k$ for $H_n(X)/\Tor H_n(X)$. Let
$q:\bar X\to X$ be the universal cover of $X$ and let $\tau$ be
the transfer. We claim that $p^ma$ can be represented by a
spherical cycle for every $a\in H_n(X)$. By the Hurewicz Theorem
every element $H_n(\bar X)$ can be represented by a spherical
cycle. Note that $p^ma=q_*(\tau_*(a))$ and the claim follows. We
attach $n+1$-cells along spherical cycles $p^ma_1,\dots, p^ma_k$
to obtain $Y$. We note $H_i(Y)=H_i(X)\in{\sC}_p$ for $i<n$. By
construction, $$H_n(Y)=\Tor
H_n(X)\oplus(\oplus_{i=1}^k\Z_{p^m})\in{\sC}_p.$$ We consider the
homology exact sequence for the pair $(Y,X)$:
$$
0\to H_{n+1}(Y)\to H_{n+1}(Y,X) @>\partial>> H_n(X)\to H_n(Y)\in{\sC}_p.
$$
the group $H_{n+1}(Y,X)$ is free and by construction its rank is
the same as that of
$H_n(X)$. Since $\partial$ is a ${\sC}_p$-epimorphism, it is
${\sC}_p$-isomorphism. Thus, $H_{n+1}=0$.
\end{pf}

PROOF OF LEMMA~\ref{block2} $(k=1)$.  As in the case $k>1$ we will
need three extra conditions (1)-(3) to run the induction. The
condition (3) remains the same. The conditions (1)-(2) are changed to
the following:

(1) {\it For every simply connected complex $K$ over the simplex
$\Delta^n$, $\nu:K\to\Delta^n$, the 1-st integral homology group
of the pull-back $\tilde K$ is isomorphic to the direct sum
$\oplus\Z_p$.}

(2) {\it for every complex $K$ over the
simplex $\Delta^n$,
$\chi:K\to\Delta^n$ the homomorphism
$f'_*:H_*(\widetilde K;\Z[\frac{1}{p}])\to H_*(K;\Z[\frac{1}{p}])$ induced by
$f'$ from the pull-back diagram
$$
\CD
\widetilde K @>\chi'>> L_n\\
@Vf'VV @Vf_nVV\\
K @>\chi>> \Delta^n.\\
\endCD
$$
is an isomorphism,
}

For $n=1$ we set $L_1=\Delta^1$ and $f_1=id_{\Delta^1}$. For $n=2$
we define $L_2$ as the mapping cone of a map $g:S^1\to S^1$ of
degree $p$. It is easy to check that all conditions are satisfied.

Assume that $f_n:L_n\to\Delta^n$ is constructed for $n\ge 2$. Let
$\Delta$ be the standard $n+1$-dimensional simplex and let
$\chi:\beta\partial\Delta\to\Delta^n$ be the characteristic map of
the barycentric subdivision. Consider
$\widetilde{\partial\Delta}$, the pull-back of the map $\chi$ and
$f_n:L_n\to\Delta^n$. Let
$f':\widetilde{\partial\Delta}\to\partial\Delta$ be the
projection. First we attach finitely many 2-cells to
$\widetilde{\partial\Delta}$ to obtain a complex $Y$ with the
abelian fundamental group. We show that $\Tor H_n(Y)\in{\sC}_p$. Indeed,
from exact sequence of the pair $(Y,\widetilde{\partial\Delta})$
follows that $\Tor H_n(Y)=\Tor H_n(\widetilde{\partial\Delta})$.
By condition (2) we have that
$H_n(\widetilde{\partial\Delta})\otimes\Z[\frac{1}{p}]\cong
H_n(\partial\Delta)\otimes\Z[\frac{1}{p}]=\Z[\frac{1}{p}]$.
Therefore $\Tor H_n(\widetilde{\partial\Delta}\in{\sC}_p$.

By condition (1) we have that
$H_1(\widetilde{\partial\Delta})=\pi_1(Y)=\oplus\Z_p$.  The group
$\pi_2(Y)$ is finitely generated since it is equal to the group
$\pi_2(\bar Y)=H_2(\bar Y)$ where $\bar Y$ is the universal cover.
We attach 3-cells killing $\pi_2(Y)$ to $Y=X_0=X_1=X_2$ to obtain $X_3$.
Similarly, the group $\pi_3(X_3)$ is finitely generated. We kill
it by attaching 4-cells and so on. We construct a chain
$X_0\subset\dots\subset X_n$ such that $\pi_1(X_n)=\oplus\Z_p$,
$\pi_i(X_n)=0$ for $2\le i<n$.
Then by Proposition~\ref{modification} $\Tor
H_n(X_n)\in{\sC}_p$. Then using Proposition~\ref{futher modif}
we attach finitely many $n+1$-cells to $X_n$ to obtain the complex
$L_{n+1}$. We define a map $f_{n+1}:L_{n+1}\to\Delta$ by sending
all new open cells in the interior of $\Delta$ by a map simplicial
with respect to some symmetric subdivision of $\Delta$.

We verify that $cd_{\Z_p}f_{n+1}=\overline{cd}_{\Z_p}f_{n+1}=1$,
$cd_{\Q}f_{n+1}=n+1$, and the conditions (1)-(3).

First we show that $\overline{cd}_{\Z_p}f_{n+1}\le 1$.  Let
$\sigma\subset\Delta^{n+1}$ be a face. If
$\sigma\subset\partial\Delta^{n+1}$, then the restriction of
$\chi$ to $\sigma$ turns it into a complex over the simplex
$\Delta^n$. Then by the induction assumption (3), the inclusion
homomorphism $$H^1((f')^{-1}(\sigma);\Z_p)\to
H^1((f')^{-1}(\partial\sigma);\Z_p)$$ is an epimorphism. If
$\sigma=\Delta^{n+1}$, the inclusion homomorphism
$H^1(L_{n+1};\Z_p)\to H^1(\widetilde{\partial\Delta};\Z_p)$ is an
epimorphism as a composition of epimorphisms
$$H^1(L_{n+1};\Z_p)\to H^1(X_n;\Z_p)\to
...\to H^1(X_3;\Z_p)\to H^1(Y;\Z_p)\to
H^1(\widetilde{\partial\Delta};\Z_p).$$ The last homomorphism in
this chain is an epimorphism since it is dual to a monomorphism
induced by the inclusion of a complex to its abelianization. All
other homomorphisms are epimorphisms by the dimensional reason.

Clearly, ${cd}_{\Z_p}f_{n+1}\ge 1$. Therefore,
$cd_{\Z_p}f_{n+1}=\overline{cd}_{\Z_p}f_{n+1}=1$.

By construction, $H_n(L_{n+1})\in{\sC}_p$. Hence the inclusion
$\widetilde{\partial\Delta^{n+1}}\subset L_{n+1}$ induces zero
homomorphism in $n$-dimensional $\Q$-homology and hence in
$\Q$-cohomology. This implies that $cd_{\Q}f_{n+1}= n+1$.

(1). Let $\chi:K\to\Delta^{n+1}$ be a light simplicial map and $K$
is simply connected. Since $n\ge 2$, the $n$-skeleton $K^{(n)}$ is
simply connected. By induction assumption
$H_1((f')^{-1}(K^{(n)}))=\oplus\Z_p$. Note that
$\widetilde{K}/((f')^{-1}(K^{(n)}))=\vee
(L_{n+1}/\widetilde{\partial\Delta})$ is the wedge of simply
connected CW complexes. From exact sequence of the pair
$(\widetilde{K},(f')^{-1}(K^{(n)}))$ it follows that
$H_1((f')^{-1}(K^{(n)}))\to H_1(\widetilde{K})$ is an epimorphism.
Hence $H_1(\widetilde{K})=\oplus\Z_p$.

(2). Let $\nu:K\to\Delta^{n+1}$ be a light simplicial map. By the
induction assumption
$$
(f'|_{...})_*:H_*((f')^{-1}(K^{(n)});\Z[\frac{1}{p}])\to
H_*(K^{(n)};\Z[\frac{1}{p}]))$$ is an isomorphism. Consider the
diagram generated by the exact sequences of pairs and the map
$f':(\widetilde K,(f')^{-1}(K^{(n)}))\to(K,K^{(n)})$.
$$ \CD H_i(\widetilde K;\Z[\frac{1}{p}]))
@>>>H_i(\widetilde K, (f')^{-1}(K^{(n)});\Z[\frac{1}{p}])) @>>>
H_{i-1}((f')^{-1}(K^{(n)});\Z[\frac{1}{p}]))\\
@Vf'_*VV @V{\psi}VV @V(f'|_{...})_*VV\\
H_i(K;\Z[\frac{1}{p}])) @>>> H_i(K,K^{(n)};\Z[\frac{1}{p}])) @>>> H_{i-1}(K^{(n)};\Z[\frac{1}{p}])).\\
\endCD
$$
By the construction
$$\xi=(f_{n+1})_*:H_*(L_{n+1},\widetilde{\partial\Delta^{n+1}};\Z[\frac{1}{p}]))\to
H_*(\Delta^{n+1},\partial\Delta^{n+1};\Z[\frac{1}{p}]))$$ is an
isomorphism. Therefore, $\psi$ is an isomorphism as the direct sum
of $\xi$. By Five Lemma $f'_*$ is an isomorphism.

(3). It suffices to show that every map $\phi:(f')^{-1}(A)\to
K(\Z_p,1)$ has an extension $\bar\phi:\tilde K\to K(\Z_p,1)$. As
in the proof of (2) we may use the induction assumption to
construct an extension $\phi':\widetilde{K^{(n)}}\to K(\Z_p,1)$ of
the map $\phi|_{(f')^{-1}(A\cap K^{(n)})}$. Since the inclusion
$\widetilde{\partial\Delta^{n+1}}\subset L_{n+1}$ induces an
epimorphism in $1$-dimensional mod $p$ cohomology, there is an
extension $\bar\phi:\tilde K\to K(\Z_p,1)$ of the map
$\phi\cup\phi':(f')^{-1}(A)\cup\widetilde{K^{(n)}}\to
K(\Z_p,1)$.\qed

\

The proof of Theorem~\ref{main 2} for $k=1$ is the same as for
$k>1$.

\end{document}